\documentclass{amsart}
\usepackage{amssymb,amsmath}

\newtheorem{theorem}{Theorem}
\newtheorem{lemma}{Lemma}
\newtheorem{conjecture}{Conjecture}

\begin{document}

\title{Partitions into a small number of part sizes}
\author{William J. Keith}
\keywords{partitions; overpartitions; divisor function; sums of squares}
\subjclass[2010]{05A17, 11P83}
\maketitle

\begin{abstract}

We study $\nu_k(n)$, the number of partitions of $n$ into $k$ part sizes, and find numerous arithmetic progressions where $\nu_2$ and $\nu_3$ take on values divisible by 2 and 4.  Expanding earlier work, we show $\nu_2(An+B) \equiv 0 \pmod{4}$ for (A,B) = (36,30), (72,42), (252,114), (196,70), and likely many other progressions for which our method should easily generalize.  Of some independent interest, we prove that the overpartition function $\overline{p}(n) \equiv 0 \pmod{16}$ in the first three progressions (the fourth is known), and thereby show that $\nu_3(An+B) \equiv 0 \pmod{2}$ in each of these progressions as well, and discuss the relationship between these congruences in more generality.  We end with open questions in this area.

\end{abstract}

\section{Introduction}

Denote the number of partitions of $n$ in which exactly $k$ sizes of part appear by $\nu_k(n)$.  For instance, $\nu_2(5) = 5$, counting $$4+1, 3+2, 3+1+1, 2+2+1, \text{ and } 2+1+1+1.$$  This easily stated function has been studied by Major P. A. MacMahon \cite{MacMahon}, George Andrews \cite{GEA1}, and more recently Tani and Bouroubi \cite{TandB}, the latter specifically interested in $\nu_2$.  The author in a recent paper \cite{Keith1} stated several theorems concerning $\nu_2$ and ventured further conjectures regarding $\nu_2$ and $\nu_3$, which it is the purpose of this paper to prove and expand.  Despite attention from these authors, results of the kind found in other areas of partition theory, such as congruences in arithmetic progressions, have not been forthcoming; here we provide several, with a proof strategy easily adaptable to future possible candidates.

Data on $\nu_k(n)$ relates to the study of overpartitions.  An overpartition of $n$ is a partition of $n$ in which the last appearance of a given size of summand is either overlined or not.  The overpartitions of 3 are \[ 3, \overline{3}, 2+1, \overline{2}+1, 2+\overline{1}, \overline{2}+ \overline{1}, 1+1+1, 1+1+\overline{1} .\]

Often attributed originally to Major MacMahon, overpartitions have seen a surge of interest in recent years since the 2004 publication of a paper by Corteel and Lovejoy \cite{CoLo}, placing them in the context of more recent work in partition theory.

Denote the number of overpartitions of $n$ by $\overline{p}(n)$.  Then it is clear that $$\overline{p}(n) = 2 \nu_1 (n) + 4 \nu_2 (n) + 8 \nu_3(n) + \dots .$$  Thus data about $\nu_i(n)$ can inform or be informed by results on overpartitions.  An example by Byungchan Kim \cite{Kim} is the theorem that $\overline{p}(n) \equiv 0 \pmod{8}$ if $n$ is neither a square nor twice a square; this is equivalent to the claim that for such numbers, $\frac{1}{2} \nu_1(n)$ and $\nu_2(n)$ are simultaneously both even or both odd.

Our main theorems include several on $\nu_2(An+B)$ mod 4 and $\nu_3$ mod 2:

\begin{theorem}\label{Nu2} $\nu_2(An+B) \equiv 0 \pmod{4}$ if $(A,B) \in \{(36,30), (72,42), (196,70), (252,114) \}$.
\end{theorem}

\begin{theorem}\label{Nu3} For each of the $(A,B)$ above, $\nu_3(An+B) \equiv 0 \pmod{2}$.\end{theorem}

In order to prove Theorem \ref{Nu3}, we need several overpartition congruences modulo 16.  The congruence for $(A,B) = (196,70)$ is already known; the others, and the generating function dissections we provide which prove them, are apparently new, although this is a field of active research.  We record these below.

\begin{theorem}\label{OverPtns} $\overline{p}(An+B) \equiv 0 \pmod{16}$ for all $(A,B)$ above.\end{theorem}

There are many other candidate progressions.  The proof techniques we give should be easily adaptable.

\phantom{.} 

\noindent \textbf{Remark:} Between versions of this article, a paper appeared by Xinhua Xiong \cite{Xiong} in which $\overline{p}(n)$ is completely determined modulo 16 by the factorization of $n$.  The overpartition identities in this paper would then also follow from Xiong's work and facts such as our candidate progressions never containing squares or sums of squares, along with observations regarding various primes appearing in their factorization.  This work, plus the method elaborated below for $\nu_2$, together give a method for obtaining many progressions for $\nu_3$.

\phantom{.}

In the next section we give much of the background information necessary to verify the results in this paper, including useful formulas from MacMahon and Andrews for $\nu_i$, and several facts concerning modular forms which are central the methodology.  The author sincerely thanks Jeremy Rouse for assistance provided on MathOverflow (\cite{MOJRouse}, \cite{MOJRouse2}) answering related questions, and a careful referee for correcting an oversight in an earlier draft.

In Section 3 we prove Theorem \ref{Nu2} by expanding Rouse' original method: we isolate as much as possible common to all such progressions using MacMahon and Andrews' results, reducing the proof in each individual case to a short catalog of necessary modular forms, and a Sturm calculation verifying a summatory congruence.  In Section 4 we prove Theorem \ref{Nu3}, proving Theorem \ref{OverPtns} in the process, and discuss challenges and possible avenues of attack in proceeding further.  The last section gives a number of open questions which we think are of general interest.

\section{Background Theorems}

Since partitions into exactly one size of part have Ferrers diagrams which are just rectangles of area $n$, $\nu_1(n)$ is just $d(n)$, the divisor function, which is perfectly understood.  If the factorization of $n$ into primes is $n = p_1^{\alpha_1} p_2^{\alpha_2} \dots$, then $$d(n) = (\alpha_1 + 1)(\alpha_2 + 1) \dots.$$  We are thus more interested in $\nu_2$ and $\nu_3$.

MacMahon and Andrews gave generating functions for $\nu_k$ and, along with Karl Dilcher independently \cite{Dilcher}, all derived the identities

\begin{equation}\label{Nu2Eq} \nu_2(n) = \frac{1}{2}\left(\sum_{k=1}^{n-1} d(k)d(n-k) - \sigma_1(n) + d(n) \right)\end{equation}

and

\begin{multline}\label{Nu3Eq} \nu_3(n) = \frac{1}{3} d(n) - \frac{1}{2} \sigma_1(n) + \frac{1}{6} \sigma_2(n) - \frac{1}{2} \sum_{k=1}^{n-1} d(k) \sigma_1(n-k) \\ + \frac{1}{2} \sum_{k=1}^{n-1} d(k) d(n-k) + \frac{1}{6} \sum_{k=1}^{n-2} \sum_{j=1}^{n-k-1} d(k) d(j) d(n-k-j)\end{multline}

\noindent where $\sigma_k(n) = \sum_{d \vert n} d^k$.  (Dilcher's identity is different in form but closely related.)

Using these ideas, the author showed in \cite{Keith1} that

\begin{theorem}\label{V2Mod2} If $n \equiv 2 \pmod{4}$, or $n$ has two or more primes appearing to odd order in its prime factorization, then $\nu_2(n) \equiv 0 \pmod{2}$.\end{theorem}

Together with Rouse, it was further shown that

\begin{theorem}\label{16Mod14} $v_2(16j+14) \equiv 0 \pmod{4}$.
\end{theorem}

Since the proof strategy for Theorem \ref{Nu2} is an expansion of this method, we sketch the proof for Theorem \ref{16Mod14} below.

One observes that for $n = 16j+14$, $\sigma_1(n) \equiv 0 \pmod{8}$ and that $d(n) \equiv d\left(\frac{n}{2}\right)^2 \pmod{8}$, so these can be removed from equation (\ref{Nu2Eq}) and it remains to show that $$\sum_{k=1}^{\frac{n-2}{2}} d(k) d(n-k) \equiv 0 \pmod{4}.$$

There are no odd terms, since $n$ is not the sum of two squares (observe quadratic residues mod 16), and therefore we wish to show that there are an even number of terms that are not multiples of 4.  The only terms that are not multiples of 4 are those in which $k$ or $n-k$ is square, and the other term is 2 mod 4, i.e. $n-k$ or $k$ respectively is $p y^2$ for $p$ a prime, with $s_p(y) \equiv 0 \pmod{2}$ where $s_p(y)$ is the power of $p$ in the prime factorization of $y$.

Thus the theorem reduces to showing that there are an even number of representations of $n$ in the form $n=x^2+py^2$ with the appropriate conditions on the prime $p$, since for each such pair $k$ will be the smaller of the two terms $x^2$ or $py^2$.  In order to analyze the parity of the number of such representations, we avail ourselves of the congruences

$$F(q) := \sum_{n=0}^\infty \sigma_1(2n+1) q^{2n+1} \equiv \sum_{n =1}^\infty q^{(2n+1)^2} \pmod{2}$$

and

$$G(q) := \frac{1}{2} \sum_{n=0}^\infty \sigma_1(8n+5) q^{8n+5} \equiv \sum_{{p \equiv 5 \pmod{8}} \atop {y \geq 1, 2 \vert s_p(y)}} q^{p y^2} \pmod{2}.$$

(When we state of functions $F(q) = \sum_{n=0}^\infty f(n)q^n$ and $G(q)=\sum_{n=0}^\infty g(n) q^n$ that $F(q) \equiv G(q) \pmod{c}$, we mean that $f(n) \equiv g(n) \pmod{c}$ for all $n$.)

With these functions, $$T(q) = F(q)G(q) + F(q^4)F(q^2)$$ \noindent has coefficients of parity equal to the number of representations we desire.

This construction is advantageous since $F(q)$ and $G(q)$ are modular forms, and thus, by the properties listed below, so is $T(q)$. Indeed we can calculate that $T(q)$ is a modular form of weight 4 for $\Gamma_0(64)$ and hence the Sturm bound is 32; a short calculation of the type described below shows that all coefficients are even, and thus the theorem is shown.

The facts in the preceding paragraph are due to the properties of modular forms.  We refer the interested reader to any textbook on modular forms for a more detailed study; we here summarize the properties we need.

\begin{itemize}
\item A modular form is said to be of weight $k$ for $\Gamma_0(N)$ or $\Gamma_1(N)$, certain subgroups of the modular group on the upper half-plane. Its \emph{level} is the minimum possible $N$.  Such a form is also of weight $k$ for any $\Gamma_0(cN)$ or $\Gamma_1(cN)$, $c \in \mathbb{N}$.
\item Modular forms of a given weight for $\Gamma_i(N)$ form vector spaces over $\mathbb{C}$.
\item The substitutions $q \rightarrow q^c$ for $c \in \mathbb{N}$ send forms of weight $k$ for $\Gamma_i(N)$ to forms of weight $k$ for $\Gamma_i(cN)$.
\item The product of two modular forms for $\Gamma_i(N)$ of weights $k$ and $\ell$ is a modular form for $\Gamma_i(N)$ of weight $k+\ell$.
\item For a form $f(q)$, if all coefficients of $q^i$ in $f$ for $i$ below the \emph{Sturm bound} are divisible by a given prime, then all coefficients of $f$ are so divisible.  This bound is $\frac{k}{12} N \prod_{p \vert N} \left(\frac{p+1}{p}\right)$ (the product is over all primes dividing $N$) for a form in $\Gamma_0(N)$ of weight $k$ and level dividing $N$, and for a form in $\Gamma_1(N)$ is increased by a factor equal to the index of the subgroup of $\Gamma_0(N)$ for which $f(q)$ is a form.
\item If $f(q) = \sum_{n=0}^\infty a(n) q^n$ is a modular form of weight $k$ and level $N$ for $\Gamma_0(N)$, then for $m \vert N$, $g(q) = \sum_{n=0}^\infty a(mn) q^n$ is also a modular form of weight $k$ and level $N$ for $\Gamma_0(N)$, and if $\chi$ is a primitive Dirichlet character mod $M$, then $g(q) = \sum_{n=0}^\infty a(n) \chi (n) q^n$ is a modular form of weight $k$ for $\Gamma_1(N M^2)$.
\end{itemize}

The last property allows us to dissect modular forms as needed for our proofs, for by selecting characters that cancel properly when the forms are added, we may construct from the form $f(q) = \sum_{n=0}^\infty a(n) q^n$ a form $g(q) = \sum_{n=0}^\infty a(An+B) q^{An+B}$ of the same weight and higher level.  The form constructed from such twists will likely lie in $\Gamma_1(N)$ for the required level, in which case the Sturm bound will be increased by a factor equal to the order of the subgroup of squares of Dirichlet characters of modulo $A$ in the group of all Dirichlet characters modulo $A$.

Our proofs will require the facts that $F(q)$ (defined above) is of weight 2 and level 4, and $H(q) := \sum_{n=0}^\infty \sigma_1(3n+1) q^{3n+1}$ is of weight 2 and level 9.

\section{Partitions into 2 sizes of part}

In \cite{Keith1} it was conjectured that $\nu_2(36n+30) \equiv 0 \pmod{4}$.  This and more is true.  We prove Theorem \ref{Nu2} by an expansion of the methodology above, executing the proof strategy in detail for $(A,B) = (36,30)$.

\begin{proof}
Set $n=36j+30$.  We again observe that $\sigma_1(n) \equiv 0 \pmod{8}$ (since $3 \vert\vert n$ and at least one prime $6\ell+5$ appears to odd order in its factorization) and that $d(n) \equiv d(\frac{n}{2})^2 \pmod{8}$, and so again it suffices to show 

\begin{equation}\label{ShortD2} \sum_{k=1}^{\frac{n-2}{2}} d(k) d(n-k) \equiv 0 \pmod{4}.
\end{equation}

By the same argument as before, we wish to show that there exist an even number of representations $n = x^2+py^2$, $s_p(y) \equiv 0 \pmod{2}$.  There are now six possible residue classes for $x$, with several possible values mod 36 of $p$ and $y$ for each.  We summarize these in the following table.

\begin{center}\begin{tabular}{|c|c|c|}
\hline $x^2 \pmod{36}$ & $py^2 \pmod{36}$ & Possible $(p,y^2) \pmod{36}$ \\
\hline 1 & 29 & $\{(29,1), (17,25) , (5,13)\}$ \\
\hline 13 & 17 & $\{(17,1), (5,25), (29,13)\}$ \\
\hline 25 & 5 & $\{(5,1) , (29,25), (17,13)\}$ \\
\hline 4 & 26 & $\{(2,13)\}$ \\
\hline 16 & 14 & $\{(2,25)\}$ \\
\hline 28 & 2 & $\{(2,1)\}$ \\
\hline
\end{tabular}\end{center}

We construct the following modular forms.  All $q$-series congruences are mod 2.

For $i = 1,25, \text{ or } 13$, with $\epsilon(i) = 1, 5, \text{ or } 7$ respectively,

$$F_{x,i}(q) := \sum_{j=0}^\infty \sigma_1(36j+i) q^{36j+i} \equiv \sum_{{j=0} \atop {j \equiv \pm \epsilon(i) \, (\text{mod }18})}^\infty q^{j^2}.$$

To illustrate for clarity, $F_{x,13} = \sum_{j=0}^\infty \sigma_1(36j+13) q^{36j+13}$, which is congruent modulo 2 to $\sum q^{j^2}$ with the sum taken over positive $j \equiv \pm 7 \pmod{18}$.

For $\ell = 4,16, \text{ or } 28$, with $\epsilon(\ell) = 2, 4, \text{ or } 8$ respectively,

$$F_{x,\ell}(q) := \sum_{j=0}^\infty \sigma_1(9j+\ell/4) \left(q^4\right)^{9j+\ell/4} \equiv \sum_{j=0}^\infty \sigma_1(36j+\ell) q^{36j+\ell} \equiv \sum_{{j=0} \atop {j \equiv \pm \epsilon(\ell) \, (\text{mod } 18)}}^\infty q^{j^2}.$$

For $m = 13, 7, \text{ or } 1$, we observe that $\sigma_1(18j+m) \equiv \sigma_1(36j+2m) \pmod{2}$.  Let $\epsilon(m) = 7, 5, \text{ or } 1$, respectively.  We construct for these $m$

$$G_{y,2m}(q) := \sum_{j=0}^\infty \sigma_1(18j+m) \left(q^2\right)^{18j+m} \equiv \sum_{j=0}^\infty \sigma_1(36j+2m) q^{36j+2m} \equiv \sum_{{j=0} \atop {j \equiv \pm \epsilon(m) \, (\text{mod } 18)}}^\infty \left(q^2\right)^{j^2}.$$

Finally, for $k = 29, 17, \text{ or } 5$, define 

\begin{align*}
G_{y,29}(q) :=& \frac{1}{2} \sum_{j=0}^\infty \sigma_1(36j+29) q^{36j+29} \equiv \sum_{{p \equiv 29 \pmod{36}} \atop {{y \equiv \pm 1 \pmod{18}} \atop {s_p(y) \, \text{even}}} } q^{p y^2} + \sum_{{p \equiv 17 \pmod{36}} \atop {{y \equiv \pm 5 \pmod{18}} \atop {s_p(y) \, \text{even}}} } q^{p y^2} + \sum_{{p \equiv 5 \pmod{36}} \atop {{y \equiv \pm 7 \pmod{18}} \atop {s_p(y) \, \text{even}}} } q^{p y^2} \\
G_{y,17}(q) :=& \frac{1}{2} \sum_{j=0}^\infty \sigma_1(36j+17) q^{36j+17} \equiv \sum_{{p \equiv 17 \pmod{36}} \atop {{y \equiv \pm 1 \pmod{18}} \atop {s_p(y) \, \text{even}}} } q^{p y^2} + \sum_{{p \equiv 5 \pmod{36}} \atop {{y \equiv \pm 5 \pmod{18}} \atop {s_p(y) \, \text{even}}} } q^{p y^2} + \sum_{{p \equiv 29 \pmod{36}} \atop {{y \equiv \pm 7 \pmod{18}} \atop {s_p(y) \, \text{even}}} } q^{p y^2} \\
G_{y,5}(q) :=& \frac{1}{2} \sum_{j=0}^\infty \sigma_1(36j+5) q^{36j+5} \equiv \sum_{{p \equiv 5 \pmod{36}} \atop {{y \equiv \pm 1 \pmod{18}} \atop {s_p(y) \, \text{even}}} } q^{p y^2} + \sum_{{p \equiv 29 \pmod{36}} \atop {{y \equiv \pm 5 \pmod{18}} \atop {s_p(y) \, \text{even}}} } q^{p y^2} + \sum_{{p \equiv 17 \pmod{36}} \atop {{y \equiv \pm 7 \pmod{18}} \atop {s_p(y) \, \text{even}}} } q^{p y^2}.
\end{align*}

These modular forms have the following weights and levels: $F_{x,1}$, $F_{x,25}$, $F_{x,13}$, $G_{y,29}$, $G_{y,17}$, and $G_{y,5}$ are all dissections of $F(q)$ by characters mod 36, and so they are all of weight 2 for $\Gamma_1(36^2 \cdot 4) = \Gamma_1(5184)$.  $F_{x,4}$, $F_{x,16}$, and $F_{x,28}$ are dissections of $H(q)$ by characters mod 9, thereafter magnified by the substitution $q\rightarrow q^4$, so they of weight 2 for $\Gamma_1(4\cdot 9^3) = \Gamma_1(2916)$.  $G_{y,26}$, $G_{y,14}$ and $G_{y,2}$ are all dissections of $H(q)$ by characters mod 18, thereafter magnified by the substitution $q \rightarrow q^2$, so they are modular forms of weight 2 for $\Gamma_1(9 \cdot 18^2 \cdot 2) = \Gamma_1(5832)$.

The product of any two of these modular forms of weight 2 for $\Gamma_1(N_1)$ and $\Gamma_1(N_2)$ is a modular form of weight 4 for $\Gamma_1(lcm(N_1,N_2))$.  For the odd $F$ and $G$, we have $F_{x,i} G_{y,j}$ of weight 4 for $\Gamma_1(5184)$.  The even cases $F_{x,2i} G_{y,2j}$  are of weight 4 for $\Gamma_1(5832)$.  Finally, the sum of all these is a modular form of weight 4 for $\Gamma_1(lcm(5184,5832))=\Gamma_1(46656)=\Gamma_1(6^6)$.

Define $S=\{1,25,13,4,16,28\}$ and set $$R(q) = \sum_{n=0}^\infty r(n) q^n = \sum_{i \in S} F_{x,i} (q) G_{y,30-i} (q).$$  Then $r(n)$ is 0 for terms other than $n=36j+30$, and for these terms is of the same parity as the number of representations of $n$ of the form sought.

If $R(q)$ were in $\Gamma_0(N)$, the Sturm bound for $R(q)$ would be $\frac{4}{12} 46656 \left(\frac{3}{2}\frac{4}{3}\right) = 31104$.  However, $R(q)$ instead lies in $M_2(\Gamma_0(6^6)) \oplus M_2(\Gamma_0(6^6), \chi) \oplus M_2(\Gamma_0(6^6),\chi^2)$ where $\chi$ is some Dirichlet character of order 3, and thus $R(q)$ is a modular form of weight 4 for a subgroup of $\Gamma_0(6^6)$ of index 3, so the bound required is 93312.

It is a straightforward calculation to construct all these forms in Mathematica or another symbolic computation package, expand the series to the Sturm bound, and check that all coefficients up to $q^{93312}$ are even.  Hence, all coefficients are even, and so the number of representations of $36j+30$ of the form required is also even.  Thus, $\sum_{k=1}^{\frac{n-2}{2}} d(k) d(n-k) \equiv 0 \pmod{4}$, and the theorem holds.

For progressions $n=Aj+B$ with $A$ even in which $\sigma_1(n) \equiv 0 \pmod{8}$ and $d(n) \equiv d\left(\frac{n}{2}\right)^2 \pmod{8}$, the sum reduces the same way to showing $$\sum_{k=1}^{\frac{n-2}{2}} d(k) d(n-k) \equiv 0 \pmod{4}.$$  If the candidate progression is among those which can never contain sums of two squares, then we again reduce the question to analyzing the parity of the number of representations of $n$ of the form $x^2+py^2$, $s_p(y) \equiv 0 \pmod{2}$, which if $Aj+B$ is a suitable progression we can analyze exactly as before.

The other progressions of the theorem can be analyzed in such a fashion.  We omit the repetitive details, noting only that for $A=72, 196, 252$ our multiples of the $\Gamma_0$ bound are 6, 21, and 9 respectively; the necessary parity checks can be easily verified by a symbolic computation package on a laptop computer. \end{proof}

\noindent \textbf{Remarks:} These are not exhaustive even of candidates of small moduli.  We note that computation has not yet suggested a candidate progression in which the conditions described do \emph{not} hold.  It would be reasonable to conjecture that the conditions are necessary:

\begin{conjecture}\label{Nu2Conj} If $\nu_2(An+B) \equiv 0 \pmod{4}$ for all $n$, then for all $n$ it holds that $\sigma_1(An+B) \equiv 0 \pmod{8}$, $d(An+B) \equiv d\left(\frac{An+B}{2}\right)^2 \pmod{8}$, and the progression $An+B$ never contains sums of two squares.
\end{conjecture}

\section{Partitions into 3 sizes of part}

We observed in the introduction that the result of Kim, that $\overline{p}(n) \equiv 0 \pmod{8}$ for $n \neq k^2, 2k^2$, is equivalent to the claim that for such $n$, $\frac{1}{2}\nu_1(n)$ and $\nu_2(n)$ are simultaneously both even or both odd.  Once we have information on $\nu_i$ for $i < k$, information about overpartitions gives us information about $\nu_k$.  We will prove Theorem \ref{Nu3} by first giving several facts about $\nu_1$ and $\nu_2$, then proving or employing congruences for $\overline{p}(An+B)$ modulo 16.

\phantom{.}

\noindent \emph{Proof of Theorem \ref{Nu3}.} Begin by observing that for $(A,B)$ one of the ordered pairs $\{(36,30), (72,42), (196,70), (252,114) \}$, $\nu_1(An+B) \equiv 0 \pmod{8}$ because at least three primes divide $An+B$ with odd exponent.  To wit, $36n+30 = 6(6n+5)$, and 5 is a quadratic nonresidue modulo 6, so 2, 3, and some additional prime divide $36n+30$ to odd order.  For $72n+42 = 6(12n+7)$, 7 is a quadratic nonresidue modulo 12; for $196n+70 = 14(14n+5)$, 5 is a quadratic nonresidue modulo 14; for $252n+114 = 6(42n+19)$, 19 is a quadratic nonresidue modulo 42.

We previously showed that $\nu_2(An+B) \equiv 0 \pmod{4}$ in each of these progressions.  Therefore, we have $$\overline{p}(An+B) \equiv 2\cdot 8 + 4 \cdot 4 + 8 \cdot \nu_3(An+B) + \dots \pmod{16}.$$

Thus, if $\overline{p}(An+B) \equiv 0 \pmod{16}$ in these progressions, it must follow that $\nu_3(An+B) \equiv 0 \pmod{2}$.  Hence we show Theorem \ref{OverPtns}.

\begin{proof}

The case $(A,B) = (196,70)$ is separate and follows immediately from Theorem 1.2 of Chen et. al in \cite{CHSZ}, which holds that $\overline{p}(7n) \equiv 0 \pmod{16}$ unless $7 \vert n$.

For the remaining cases, in which $36 \vert A$, we employ several identities from \cite{XiaYao}, beginning with congruence (4.29) of that paper.  For compactness, we employ the notation $$f_i = \prod_{k=1}^\infty (1-q^{ik}).$$  It holds for $\ell > 0$ that ${f_i}^{2^\ell} \equiv {f_{2i}}^{2^{\ell-1}} \pmod{2^\ell}$, and more generally for $\ell > k \geq 0$ that $$2^k {f_i}^{2^{\ell-k}} \equiv 2^k {f_{2i}}^{2^{\ell-k-1}} \pmod{2^{\ell-k}}.$$

We will require two lemmas from \cite{CHSZ}.  First is the 3-dissection of the overpartition function:

\begin{lemma}\label{OPmod3} $\frac{f_2}{{f_1}^2} = \frac{f_6^4 f_9^6}{f_3^8f_{18}^3} + 2q\frac{f_6^3f_9^3}{f_3^7}+4q^2\frac{f_6^2f_{18}^3}{f_3^6}$
\end{lemma}

Next is the 2-dissection of another quotient:

\begin{lemma}\label{ThreeEven} $\frac{f_3^3}{f_1} = \frac{f_4^3f_6^2}{f_2^2f_{12}} + q \frac{f_{12}^3}{f_4}$
\end{lemma}

Now, from equation (4.29) in \cite{CHSZ}, we extract the even terms and reduce the congruence to one modulo 16 to obtain

$$\overline{p}(6n) q^n \equiv \frac{f_2^4 f_{12}^{15}}{f_1^8 f_{6}^6 f_{24}^6} + 12q^3 \frac{f_{12}^3 f_{24}^2}{f_6^2} \pmod{16}.$$

We wish to extract from this identity those terms in which $n \equiv 5 \pmod{6}$.  In the second summand, all $q^n$ have $n \equiv 0 \pmod{3}$, so we discard these.  Take $\frac{f_2^4}{f_1^8}$ and employ Lemma \ref{OPmod3} to obtain

\begin{equation}\label{Fourth}\sum_{n=0}^\infty \overline{p}(6n)q^n \equiv \frac{f_{12}^{15}}{f_6^6 f_{24}^6} \left( \frac{f_6^4 f_9^6}{f_3^8f_{18}^3} + 2q\frac{f_6^3f_9^3}{f_3^7}+4q^2\frac{f_6^2f_{18}^3}{f_3^6} \right)^4 + \cdots \pmod{16}\end{equation}

\noindent where the elided terms do not have power $n \equiv 5 \pmod{6}$.  Now expand the fourth power, disregarding all terms with an integer coefficient divisible by 16, and extract all those terms in which $n \equiv 2 \pmod{3}$:

$$\sum_{n=0}^\infty \overline{p}(18n+12)q^{3n+2} \equiv \frac{f_{12}^{15}}{f_6^6 f_{24}^6} \left( 24q^2 \frac{f_6^7 f_9^{18}}{f_3^{30} f_{18}^6}  \right) \pmod{16}.$$

But since $24\frac{1}{f_3^{30}} \equiv 24 \frac{1}{f_6^{15}} \pmod{16}$, all powers in this congruence with odd coefficient are in fact congruent to 2 mod 6, and so no terms congruent to 5 mod 6 appear with nonzero coefficient modulo 16.  The theorem is proved for $(A,B) = (36,30)$.

To prove the case $(A,B) = (72,42)$, we start from equation \ref{Fourth} and extract terms $7 \pmod{12}$.  Begin with terms congruent to $1 \pmod{3}$:

$$\sum_{n=0}^\infty \overline{p}(18n+6)q^{3n+1} \equiv \frac{f_{12}^{15}}{f_6^6 f_{24}^6} \left( 8 q \frac{f_6^{15} f_9^{21}}{f_3^{31} f_{18}^9}  \right) \equiv 8q \frac{f_9^3}{f_3} \pmod{16}.$$

We now employ Lemma \ref{ThreeEven} to obtain those terms that are 1 mod 6:

$$ \sum_{n=0}^\infty \overline{p}(36n+6) q^{6n+1} \equiv 8q \frac{f_{12}^3 f_{18}^2}{f_6^2 f_{36}} \equiv 8q f_{24} \pmod{16}.$$

Hence there are no powers $q^j$ with $j \equiv 7 \pmod{12}$ that have coefficients nonzero modulo 16.

Finally, to show the case $(A,B) = (252,114)$, we take the congruence above: 
$$ \sum_{n=0}^\infty \overline{p}(6(6n+1)) q^{n} \equiv 8 f_{4} \pmod{16}.$$
\noindent and extract terms where $n \equiv 3 \pmod{7}$.  But by the Pentagonal Number Theorem $f_4 = \sum_{n=-\infty}^\infty (-1)^n q^{2n(3n+1)},$ and for integer argument $2n(3n+1)$ only takes on residues 0, 1, 4, or 6 modulo 7.  Thus, in the progression $(A,B) = (252,114)$ no coefficients are nonzero mod 16, and Theorem \ref{OverPtns} is proven. \end{proof}

Since Theorem \ref{OverPtns} holds, and the necessary conditions on $\nu_1$ and $\nu_2$ are fulfilled, it follows that $\nu_3 \equiv 0 \pmod{2}$ in the progressions studied.  Thus Theorem \ref{Nu3} is proved. \hfill $\Box$

\section{Open Questions}

One could possibly use information about $\nu_k$ to prove statements about overpartitions, at least modulo powers of 2.  In order to do so one would need to analyze $\nu_k$ without invoking overpartition congruences, analyzing the parity of the terms in the generating function as we did for Theorem \ref{Nu2}.  At present we require information about $\nu_1$ and $\nu_2$ to obtain information on $\nu_3$.  While it is conceivable that $\nu_3$ might possess arithmetic progressions in which all values are even without the same holding true for higher powers for $\nu_2$ and $\nu_1$, we believe there is good reason to think that this is not the case, and formally conjecture:

\begin{conjecture}\label{Nu3Conj} If $\nu_3(An+B) \equiv 0 \pmod{2}$ for some arithmetic progression, it also holds that $\nu_2(An+B) \equiv 0 \pmod{4}$.
\end{conjecture}

Why might this conjecture hold?  Suppose one wishes to show $\nu_3(36j+30) \equiv 0 \pmod{2}$ by analyzing the parity of the terms in equation (\ref{Nu3Eq}).  Suppose we have already shown for $n \equiv 30 \pmod{36}$ that $\nu_2(n) \equiv 0 \pmod{4}$, and we know $d(n) \equiv 0 \pmod{8}$, $\sigma_1(n) \equiv 0 \pmod{8}$, and $\sum_{k=1}^{n-1} d(k) d(n-k) \equiv 0 \pmod{8}$.  (In any other arithmetic progression, if any three of these are true, all four are, because we may subtract the other terms in equation (\ref{Nu2Eq}) from $\nu_2(n)$.)

We may then subtract these terms from equation (\ref{Nu3Eq}) for $\nu_3(n)$ to obtain

\begin{multline*} \nu_3(n) = \frac{1}{3} d(n) - \frac{1}{2} \sigma_1(n) + \frac{1}{6} \sigma_2(n) - \frac{1}{2} \sum_{k=1}^{n-1} d(k) \sigma_1(n-k) \\ + \frac{1}{2} \sum_{k=1}^{n-1} d(k) d(n-k) + \frac{1}{6} \sum_{k=1}^{n-2} \sum_{\ell=1}^{n-k-1} d(k) d(\ell) d(n-k-\ell) \\
\equiv -\frac{1}{6}d(n) +\frac{1}{6} \sigma_2(n) - \frac{1}{2} \sum_{k=1}^{n-1} d(k) \sigma_1(n-k) \\ + \frac{1}{6} \sum_{k=1}^{n-2} \sum_{\ell=1}^{n-k-1} d(k) d(\ell) d(n-k-\ell) \pmod{2}.
\end{multline*}

It is not difficult to show that $d(36j+30) \equiv -\sigma_2(36j+30) \pmod{12}$ (both functions being multiplicative, one simply observes the values mod 6 of each factor) and hence we can write

\begin{multline*} \nu_3(n) \equiv -\frac{1}{3}d(n) - \frac{1}{2} \sum_{k=1}^{n-1} d(k) \sigma_1(n-k) \\ + \frac{1}{6} \sum_{k=1}^{n-2} \sum_{\ell=1}^{n-k-1} d(k) d(\ell) d(n-k-\ell) \pmod{2}.
\end{multline*}

We now note that many terms in the final sum can be cast out modulo 2.  If exactly one of $k$, $\ell$, or $n-k-\ell$ is a nonsquare and the other two terms are not equal, we can group the six permutations of the three entries, the product of which are even, and discard them.  If exactly one entry is a square, we can again do so -- we may have only three permutations, but the product is divisible by 4.

If all three are non-squares, the only term we cannot permute and cast out is when $k = \ell = 12j+10$, which may not have $d(12j+10) \equiv 0 \pmod{3}$.  But $d(12j+10) = \frac{1}{2} d(36j+30)$.  If we add one-sixth of the cube of this to $-\frac{1}{3}d(n)$, we obtain $\frac{1}{3} \cdot \frac{1}{16} d(n)(d(n)-4)(d(n)+4)$.  In the latter product one term is divisible by 3 and since $d(n) \equiv 0 \pmod{8}$, the whole is an even integer.

If all three are squares, then they cannot be the same square (10 is not a quadratic residue mod 12) and thus they have 3 or 6 permutations; however, we may not be able to cast out such terms.  For instance, $30=25+4+1$ and the six permutations thereof, and this is the only such representation of 30.  We are also left with terms in which exactly one entry is a non-square and the other two are equal squares.  We may multiply by 3 and take the representative of these in which $k = \ell$ are the squares.  Thus, we end up interested in representations of $n$ by three squares, or twice a square and a non-square.  (It is interesting that overpartition identities so often relate to identities concerning sums of squares.)

Now observe that in $\frac{1}{2} \sum_{k=1}^{n-1} d(k) \sigma_1(n-k)$, $\sigma_1(n-k) \equiv d(n-k) \pmod{2}$ unless $n-k$ is twice a square, and we already know that $\frac{1}{2} \sum_{k=1}^{n-1} d(k) d(n-k) \equiv 0 \pmod{4}$.  Thus we can reduce the sought identity to 

\begin{multline} \nu_3(n) \equiv -\frac{1}{2} \sum_{k=1}^{n-1} d(k)\sigma_1(n-k) \\ + \sum_{{j+k+\ell=n} \atop {0 < j < k < \ell \text{ distinct squares}}} d(j)d(k)d(\ell) + \frac{1}{2} \sum_{k=1}^{\lfloor \sqrt{(n-1)/2} \rfloor} d(k^2)^2 d(n-2k^2) \pmod{2}.\end{multline}

We know this is congruent to 0, by the previous work; the search for a direct proof seems like a natural question of interest.

Finally, in addition to the conjectures stated previously, a number of open questions present themselves.

\begin{enumerate}
\item Treat candidate progressions in a more unified fashion, probably via the theory of eigenforms.  Can we show the existence of an infinite class of $(A,B)$ for which $\nu_2(An+B) \equiv 0 \pmod{4}$, and/or $\nu_3(An+B) \equiv 0 \pmod{2}$?
\item Numerical experimentation to date has found no progressions $An+B$ in which $\nu_2(An+B) \equiv 0 \pmod{N}$ for any $N$ other than 2 or 4; and none for $\nu(3)$ other than $N=2$.  If different moduli occur, they may have large progression modulus $A$.  Do these occur, and if so, how can they be efficiently found, or, are they forbidden?
\item Experimentation has yielded no progressions with nontrivial modulus for $\nu_k$ with $k > 3$.  It is plausible that these never occur, since from the formulas in Andrews \cite{GEA1} these values involve sums concerning $d(k) \sigma_2(n-k)$, and $\sigma_2(j)$ is not part of the same framework of modular forms and their symmetries as $\sigma_1$.  (When it appeared in $\nu_3(36j+30)$ it was a single term which cancelled with $d(n)$.)  Again, can these occur, and if so where, or are they forbidden?
\item Elaborate on the relationships between $\nu_1$ and $\nu_2$, and between $\nu_2$ and $\nu_3$.  State conditions necessary and/or sufficient for simultaneous congruences.
\item Complete the combinatorial proof for $\nu_3(36j+30)$ and generalize to other progressions.
\end{enumerate}

\section{Acknowledgements}

The author sincerely thanks Jeremy Rouse for assistance provided on MathOverflow \cite{MOJRouse, MOJRouse2}, and a careful anonymous referee for correcting an oversight in an earlier draft.

\end{document}